\documentclass[a4paper,11pt]{article}
\usepackage{geometry}
\usepackage[utf8]{inputenc}
\usepackage[english]{babel}

\catcode`\"=13 
\def"#1"{``#1''} 

\usepackage{authblk}

\usepackage{dsfont} 
\usepackage{amsmath,amssymb,amsthm}
\usepackage[justification=centering]{caption}
\usepackage[colorlinks=true, citecolor=red, linkcolor=blue]{hyperref}
\usepackage[nameinlink, capitalize]{cleveref}

\usepackage[colorinlistoftodos]{todonotes}
\usepackage{varwidth}

\usepackage{color}
\usepackage{float}
\usepackage[super]{nth}
\usepackage{booktabs}
\usepackage{rotating}
\usepackage{verbatim}
\usepackage{multirow}
\usepackage{physics}
\usepackage{mathtools}
\usepackage{subcaption}
\usepackage{bm}
\usepackage[toc,page]{appendix}
\usepackage[shortlabels]{enumitem}
\usepackage{nomencl}

\newcounter{mnotecount}[section]

\renewcommand{\themnotecount}{\thesection.\arabic{mnotecount}}

\newcommand{\mnote}[1]
{\protect{\stepcounter{mnotecount}}$^{\mbox{\footnotesize
$
\bullet$\themnotecount}}$ \marginpar{
\raggedright\em
$\!\!\!\!\!\!\,\bullet$\themnotecount: #1} }

\newcommand{\jm}[1]{{\color{purple}\mnote{{\color{purple}
#1} }}}

\def\eps {{\epsilon}}

\theoremstyle{plain}
\newtheorem{theorem}{Theorem}
\newtheorem{corollary}{Corollary}
\newtheorem{lemma}{Lemma}
\newtheorem{proposition}{Proposition}

\theoremstyle{definition}

\newtheorem{example}[theorem]{Example}
\newtheorem{definition}{Definition}
\newtheorem{remark}{Remark}

\begin{document}

\title{Velocity Averaging Lemmas: Classical,\\ Quantum and Semi-Classical}

\author[a]{François Golse}
\author[b]{Norbert J. Mauser}
\author[a,b,c]{Jakob Möller}
\affil[a]{École polytechnique \& IP Paris, CMLS\\

F-91128 Palaiseau Cedex, France\\

\texttt{francois.golse@polytechnique.edu}\vspace{5pt}}
\affil[b]{Research platform MMM c/o Fak. Math. Univ. Wien\\

Oskar Morgenstern Platz 1, A-1090 Wien, Austria\\

\texttt{norbert.mauser@univie.ac.at} \vspace{5pt}}
\affil[c]{Wolfgang Pauli Inst. c/o Fak. Math. Univ. Wien\\

Oskar Morgenstern Platz 1, A-1090 Wien, Austria\\

\texttt{jakob.moeller@univie.ac.at}\vspace{5pt}}
%

\date{\today}

\maketitle
\begin{abstract}

Averaging lemmas were introduced as a tool of the mathematical analysis of kinetic equations, i.e. PDEs for functions in phase space $(x,v)$ containing a transport ("advection") term. By integrating over $v$ in velocity space $\mathbb{R}_v^d$ (velocity averaging), one gains regularity for the density in position space $\mathbb{R}_x^d$.
The concept was invented independently by V.I. Agoshkov and by F. Golse, B. Perthame, R. Sentis and P.-L. Lions, and successfully applied to the analysis of Vlasov or Boltzmann equations in "classical kinetic theory".
In "quantum kinetic theory", the Schrödinger equation for the complex-valued "wave function" in the physical space is converted into the Wigner equation for the real-valued Wigner function in phase space (which can take negative values). The Wigner ("Quantum Vlasov") equation contains the transport term of classical kinetic equations plus a pseudo-differential operator containing the potential.
We give answers to the long standing question of whether and to which extent averaging lemmas apply to the "quantum" case of the Wigner equation.
The hard part are the "semi-classical" averaging lemmas, where one considers the asymptotics of vanishing Planck constant towards the non-negative Wigner measure. In that context "pure vs. mixed states" play a crucial role, as well as the connection between the Schrödinger equation and Quantum Hydrodynamics (QHD).
We present the results for the classical and quantum cases and sketch the "semi-classical" case which is worked out in full detail in a follow-up article.

%
%

\bigskip


\end{abstract}


\section{Introduction}


\subsection*{a) Classical velocity averaging}

"Kinetic equations" like the Vlasov or the Boltzmann equation – describing e.g. (rarefied) gases and plasmas – are time-dependent PDE for a phase space density $f\equiv f(t,x,v)\ge 0$ (also called the \emph{distribution function}), with time $t \in \mathbb R$ and $(x,v)$ in the "phase space" $\mathbb R^d_x \times \mathbb R^d_v$. Here, $x$ is the position variable, and $v$ is the velocity (resp. momentum or a more general "kinetic variable", also called $\xi$ or $p$).

Kinetic PDEs typically contain the free transport operator $\partial_t + v \cdot \nabla_x$ (also called "advection operator") which is of a hyperbolic nature.
They can be nonlinear, as in the case of the Vlasov-Poisson and Boltzmann equations. In any case the (main) unknown function in kinetic theory is the phase space density $f\equiv f(t,x,v) \ge 0$ which is the probability distribution of the "particles" at position $x$ with velocity $v$ at time $t$. 

In general, kinetic equations can be written in the form 
\begin{equation}
(\partial_t+v\cdot\nabla_x)f(t,x,v) = S[f](t,x,f)
\label{eq: general kinetic equation}
\end{equation}
where $S[f]$ is a functional of the unknown $f$, such as the collision integral in the case of the Boltzmann or the Landau equations, or a term of the form $S[f]=\nabla_v\cdot (fF[f])$, 
where $F[f]$ is a mean-field, self-consistent acceleration created by the unknown distribution $f$ itself, such as gravity in the case of the Vlasov-Poisson equation used in cosmology, or the electromagnetic force in the case of the (relativistic) Vlasov-Maxwell equation used in plasma physics.

Note that in the relativistic case the velocity $v$ is bounded by the speed of light in vacuum (we scale $c=1$) as a function $v(\xi) = \gamma(\xi) \xi$ of the unbounded momentum $\xi$ involving the "gamma factor" $ \gamma(\xi) =1/\sqrt{1+|\xi|^2}$. Therefore relativistic velocity averaging requires a slightly different setting.\\

The heuristic principle is that by "averaging in velocity", i.e. by integrating $f$ in $v$, one gains regularity in the time and position variables for the average. More precisely, if $S[f]\in L^2_{loc}(\mathbb R_t\times\mathbb R^d_x\times\mathbb R^d_v)$ 
and if $\tau+v\cdot\xi\not=0$, 
then $(t,x)\mapsto f(t,x,v)$ is microlocally in the Sobolev space $H^1$ in some conical neighborhood of $(\tau,\xi)$
around each $(t,x)\in\mathbb R\times\mathbb R^d_x$. 
On the other hand, if $f$ itself belongs to $L^2_{loc}(\mathbb R\times\mathbb R^d_x\times\mathbb R^d_v)$, 
the set $\{v\in \mathbb{R}^d \colon \tau+v\cdot\xi=0\}$ is an affine plane in $\mathbb R^d$ and has Lebesgue measure zero. 
Therefore, one can expect that there exists $s>0$ such that
\begin{equation}
\int_{|v|\le R}f(t,x,v)\dd v\in H^s_{loc}(\mathbb R_t\times\mathbb R^d_x)
\label{eq: gain in regularity}
\end{equation}
(where $H^s$ designates the Sobolev space of functions with derivatives of order $\le s$ in $L^2$) for all $R>0$.

\subsection*{b) Quantum velocity averaging\footnote{The term ``quantum averaging lemma'' has recently been coined by Han-Kwan and Rousset to designate a different, but related regularization property for the Hartree equation with Dirac
potential \cite{HanKwanRousset}.}}

By using the Wigner transform \cite{wigner1997quantum, gerard1997homogenization} of the Schr\"odinger wave function $\psi(t,x)$ (or more generally the density operator of the quantum state) one obtains the "Wigner function" $w_{\hbar}(t,x,v)$, which is the quantum mechanical analogue of the classical phase space density function $f(t,x,v)$. Here $\hbar$ is the (scaled) Planck constant, and $\hbar \to 0$ stands for the semi-classical limit.

The Wigner function $w_{\hbar}$ is a "pseudo-  phase-space density" since, on the one hand, its moments in $v$ yield the correct macroscopic quantities like the position and current densities, but, on the other hand, it assumes also negative values in general. This expresses the Heisenberg uncertainty principle that the "conjugate variables" $x$ and $v$ cannot be measured simultaneously with arbitrary precision. In fact, the Wigner function $w_{\hbar}[\psi]$ of a pure state $\psi$ is non-negative if and only if $\psi\in L^2$ corresponds to a Gaussian state (i.e. the exponential of a quadratic polynomial), as was shown by R. Hudson \cite{hudson1974wigner}. Note that Gaussian states are precisely the states for which the uncertainty is minimal (i.e. equality holds in the Heisenberg inequality). 

The Wigner function $w_{\hbar}$ obeys the "Wigner equation" (also "Quantum Vlasov equation" when coupled to a self-consistent potential), which is the quantum analogue 
to classical kinetic equations like \eqref{eq: general kinetic equation} for the phase space function. 
The Wigner transform of the Schr\"odinger equation (or the von Neumann equation) with a potential $V(t,x)$ yields the Wigner equation for the Wigner function $w_{\hbar}$
\begin{equation}
     (\partial_t + v \cdot \nabla_x) w_{\hbar}(t,x,v) = (\theta[V]w_{\hbar})(t,x,v),
     \label{eq: wigner equation}
\end{equation}
where $\theta[V]$ is a pseudodifferential operator (with a symbol involving finite differences of $V$) that acts on $w_{\hbar}$. This term can be seen as a further example of the functional $S$ in \eqref{eq: general kinetic equation} above, with the caveat that $w_{\hbar}$ in $S[w_{\hbar}](t,x,w_{\hbar})$ can now take negative values.
In the semi-classical limit $\hbar \rightarrow 0$, the Wigner equation converges to the Vlasov equation
\begin{equation} \label{Vlasov_equation}
     (\partial_t+ v \cdot \nabla_x) f  = \nabla_x V \cdot \nabla_{v} f,
\end{equation}
where $f(t,x,v) \geq 0$ is the nonnegative Wigner measure of $w_{\hbar}$ \cite{lions1993mesures} and $S[f] = \nabla_x V \cdot \nabla_{v} f$. 
Thus the classical kinetic structure and the classical velocity averaging are recovered in the semi-classical limit.

We present results explaining to which extent we can obtain "quantum averaging lemmas" for the quantum case with $\hbar$ fixed (i.e. up to the semi-classical limit).

\subsection*{c) Semi-classical velocity averaging}

The key question is what happens in the limit $\hbar \to 0$ and if one is able to obtain some gain in the regularity of averages in $v$ of $w_\hbar$in $t,x$ which is uniform in $\hbar$. 
Such semi-classical velocity averaging results should provide valuable information on observable densities, i.e. on moments in $v$ of the Wigner function $w_{\hbar}(t,x,v)$, and should be a valuable tool for semi-classical limits.

The first results in this direction are given in \cite{golse2025velocity}. One quickly notices that the uniformity in $\hbar$ puts considerable constraints on the class of admissible Wigner functions $w_{\hbar}$, where a distinction has to be made between pure 
and mixed states, expressed in terms of the density operator $R$.

This major difference between pure and mixed states is well known \cite{lions1993mesures, markowich1993classical} from the global-in-time semi-classical limit of the Schrödinger-Poisson equation to the Vlasov-Poisson equation using Wigner transforms. 
The case of $d\geq 2$ space dimensions is known only for  mixed states with strong assumptions on the (initial) occupation probabilities (that have to depend on the semiclassical parameter $\hbar$ !). The pure state case, addressed in \cite{zhang2002limit}
could only be treated in space dimension $d=1$ using non-unique measure-valued solutions of the Vlasov-Poisson equation. See however Proposition 2.4 for the case of dimension $d=3$ in the more recent reference \cite{FGTPaulCPAM2022} based on 
Serfaty's remarkable inequality for the repulsive Coulomb potential \cite{SerfatyDuke}.

We present a \emph{semi-classical averaging lemma} in  Theorem \ref{thm:VAWignerL2 semi-classical}, where we deal with a special class of mixed states for density operators bounded in Hilbert-Schmidt norm uniformly in the Planck constant. The case of pure states has to be treated differently and the key observation is a characterization of the Wigner transforms of pure states (Lemma \ref{L-CharRk1}).
Based on this characterization, it can be shown that sequences of pure states with density functions and current densities converging strongly in $L^2_{loc}$ and kinetic energy densities converging weakly in $L^2_{loc}$ have monokinetic Wigner measures (Proposition \ref{P-Monokinetic}). This concentration effect on the momentum variable in the Wigner measure explains why one cannot expect that velocity averaging can be applied to families of pure states in the semi-classical setting. \\

Finally, section \ref{sec:DerMadelung} explains how our characterization of the Wigner transforms of pure states in Lemma \ref{L-CharRk1} can be used to arrive at a quick derivation of the equations of Quantum Hydrodynamics (QHD) (due to Madelung) associated to pure states. This approach to the QHD equations is 
used to understand the physical meaning of the condition in Proposition \ref{P-Monokinetic} leading to monokinetic Wigner measures in the classical limit.


\section{Classical averaging lemmas}

In the \emph{classical} regime (i.e. "classical Newtonian physics", where the Planck constant $\hbar = 0$), averaging lemmas were first  studied independently in \cite{agoshkov1984} and \cite{bgp1985} in the $L^2$ case by studying the Fourier transform of $\partial_t +v\cdot \nabla$ and more systematically in \cite{golse1988regularity}. The subsequent literature is rich: the $L^p$ case was studied e.g. in \cite{diperna1991Lp, bezard1994regularite, arsenio2019maximal}. Velocity averaging fails in $L^1$ in general, due to possible concentrations in the $v$ variable, which offset the benefits of averaging \cite{golse1988regularity} (see Example 1 on p. 123 in \cite{golse1988regularity}). 

The basic velocity averaging theorem of interest in connection with the Wigner equation is an extension due to R. DiPerna and P.-L. Lions \cite{diperna1989global} of the main result 
in \cite{golse1988regularity}, which is particularly well-suited to handle Vlasov-type equations --- and was used to prove the global existence of weak solutions to the Cauchy problem for the Vlasov-Maxwell
system, for all square-integrable initial data with finite mass and energy.

\begin{theorem}\label{thm: diperna lions}
Let $n\ge 0$ and $f \in L^2(\mathbb{R}_t\times \mathbb{R}^d_x \times \mathbb{R}^d_{v})$ satisfy 
\begin{equation}
	(\partial_t +v \cdot \nabla_x)f = S[f] \quad \text{in }\mathcal{D}'(\mathbb{R}_t\times \mathbb{R}^d_x \times \mathbb{R}^d_{v})
    \label{eq:kinetic equation dipernalions}
\end{equation}
where $S[f] \in L^2(\mathbb{R}_t \times \mathbb{R}^d_x, H^{-n}(\mathbb{R}^d_{v}))$.  
Then for each $\psi\in\mathcal S(\mathbb{R}^d_v)$,
 \begin{equation}
 	\rho[f] := \int_{\mathbb R^d}f(\cdot,\cdot,v) \psi(v) \dd v \in H^{s}(\mathbb{R}_t\times \mathbb{R}^d_x), \quad  s=\frac{1}{2(n+1)}
\end{equation}
\end{theorem}
Note that in \cite{diperna1989global}, the RHS of \eqref{eq:kinetic equation dipernalions} is denoted by $g$, while we denote it by $S[f]$ in order to keep our notation consistent.

The relativistic version of this theorem is due to G. Rein \cite{rein2004global} and replaces the unbounded velocity $v$ by the relativistic velocity $v(\xi) = \xi/\sqrt{1+|\xi|^2}$ which is bounded by the speed of light $c$ (which is scaled to $c=1$).

\begin{theorem}\label{thm: rein}
Let $R>0$ and $f,g,h \in L^2(\mathbb{R}_t\times \mathbb{R}^d_x \times B_R)$ satisfy 
\begin{equation}
	\left(\partial_t +\frac{\xi}{\sqrt{1+|\xi|^2}} \cdot \nabla_x\right)f = g + \mathrm{div}_{\xi} h\quad \text{in }\mathcal{D}'(\mathbb{R}_t\times \mathbb{R}^d_x \times B_R)
    \label{eq:kinetic equation rein}
\end{equation} 
Then for each $\psi\in\mathcal S(B_R)$,
 \begin{equation}
 	\rho[f] := \int_{\mathbb R^d}f(\cdot,\cdot,\xi) \psi(\xi) \dd \xi \in H^{1/4}(\mathbb{R}_t\times \mathbb{R}^d_x).
\end{equation}
\end{theorem}


\section{Quantum averaging lemmas}

We now state in mathematical detail how to adapt the classical case to quantum physics, where the nonnegative distribution function $f(t,x,v)$ is replaced by the Wigner function $w_{\hbar}(t,x,v)$ (or the Wigner matrix $W_{\hbar}(t,x,v)$ in the case of vector-valued equations).

In the key quantum mechanical concept of "pure states" vs. "mixed states" the simple wave function is replaced by the more general "density matrix" (i.e. the integral kernel of the "density operator"), which doubles the number of independent variables, but which has the mathematical advantage that one can impose conditions on the "occupation probabilities" $\lambda_j$ of the initial mixed state, which have to depend on $\hbar$ to ensure uniform bounds in $\hbar$. This is essential for semi-classical limits and semi-classical averaging lemmas.

\subsection{Preliminaries}\label{sec:Prelim}

In quantum statistical mechanics the "mixed state" of a system, corresponding to infinitely many state vectors, is described by a \emph{density operator} $R$, a self-adjoint, nonnegative operator on $L^2(\mathbb{R}_x^d)$ such that $\tr_{L^2} R = 1$) . \\
The special case of a \emph{pure state} that can be described by a single wave-function $\psi \in L^2(\mathbb{R}_x^d)$ corresponds to a rank-one density operator $R=\ket{\psi}\bra{\psi}$. Here we use the bra-ket notation: $\ket{\psi}$ stands for the vector
$\psi\in L^2(\mathbb{R}_x^d)$, i.e. the square-integrable function $x\mapsto\psi(x)$ defined for a.e. $x\in\mathbb R^d$, while $\bra{\psi}$ stands for the covector defined by $\psi$, i.e. the continuous linear functional 
\[
\phi\mapsto\int_{\mathbb R^d}\overline{\psi(x)}\phi(x)\dd x
\]
defined on $L^2(\mathbb{R}_x^d)$.\\
A general density operator $R$ can be written as an integral operator and is usually represented by its kernel $R(X,Y)\in L^2(\mathbb{R}_X^d\times \mathbb{R}_Y^d)$, called the \emph{density matrix}, abusively denoted by the same letter. The operator $R$ 
has a sequence of real eigenvalues $(\lambda_j)_{j\ge 1}$ and a Hilbert basis of eigenfunctions $\{\psi_j\,:\,j\ge 1\}$ of $L^2(\mathbb{R}^d_x)$ such that
\begin{equation}\label{SpecDecR}
    R(X,Y) = \sum_{j\ge 1}\lambda_j \psi_j(X) \overline{\psi_j(Y)}, \qquad \lambda_j\ge 0, \qquad \sum_{j\ge 1}\lambda_j=\tr_{L^2(\mathbb R^d)}R=1.
\end{equation}
We can define the \emph{density function} $\rho\in L^1(\mathbb{R}^d)$ (on account of \eqref{SpecDecR}) and the trace of $R$ as 
\begin{equation}\label{DefDensityFunc}
    \rho(X) := R(X,X), \quad \tr_{L^2(\mathbb R^d)} R = \int_{\mathbb{R}^d} \rho(X) \dd X
\end{equation}
The quantum mechanical analogue of the Liouville equation for the phase-space density $f$ is the von Neumann equation for the density operator $R$,
\begin{equation}\label{eq:von Neumann} 
\left\{
\begin{aligned}
{}&i\hbar \partial_t R(t) = [\mathcal H,R(t)], \quad \mathcal H:=-\tfrac{\hbar^2}{2m}\Delta+V,
\\
&R(0) = R_{0}.
\end{aligned}
\right.
\end{equation}
where $[\mathcal H,R(t)]$ is the commutator of $R$ with the Hamilton \textit{operator} $\mathcal H$ (with $V$ as the usual multiplication operator), which is such that $\tfrac{i}\hbar[\mathcal H,R(t)]$ converges to the Poisson bracket of the Hamilton \textit{function}
and the classical density $f$ in the semiclassical limit.
Now consider the density matrix $\tilde{R}$ in the "center of mass" and "relative" variables $x:= \tfrac{1}{2}(X+Y)$ and $y:=\tfrac{1}{\hbar}(X-Y)$:
\begin{equation}\label{eq:weyl variable density matrix}
        \tilde{R}(t,x,y):= R(t,X,Y)=R(t,x+\tfrac\hbar{2}y,x-\tfrac\hbar{2}y).    
\end{equation}
The \emph{Wigner transform} $w_{\hbar}$ of $R$ is defined as the Fourier transform of the density matrix that turns the relative position variable $y$  into a dual "velocity variable" $v$:
\begin{equation}
    w_{\hbar}[R](t,x,v):=\tfrac1{(2\pi)^d}\mathcal{F}_{y\rightarrow v} [\tilde{R}](t,x,v)=\tfrac{1}{(2\pi)^d} \int_{\mathbb{R}_y^d} e^{-iv\cdot y} \tilde{R}(t,x,y)  \dd y.
\end{equation}
Note that we explicitely mark the $\hbar$ dependence of $w_{\hbar}$ by a subscript $\hbar$ but omit denoting the $\hbar$ dependence of $R$ for notational simplicity.

The zeroth order, resp. first order moments of the Wigner function $w_{\hbar}$ in the $v$ variable are the \emph{particle density}, resp. the vector-valued \emph{current density}
\begin{equation}
\label{eq:restriction density}
    \rho_{\hbar}(t,x) := \int_{\mathbb{R}_v^d}w_{\hbar}(t,x,v)\dd v, \quad J_{\hbar}(t,x) := \tfrac1m\int_{\mathbb{R}_v^d}v w_{\hbar}(t,x,v)\dd v. 
\end{equation}
The Wigner function $w_{\hbar}:=(2\pi)^{-d}\mathcal F_{y\to v}\tilde R_\hbar$
solves the \emph{Wigner equation}
\begin{align}
    {}&\partial_tw_{\hbar} +\tfrac1mv \cdot \nabla_x w_{\hbar} - \theta[V]w_{\hbar} = 0,\label{eq:wigner equation} 
    \\
    &w_{\hbar}(0,x,v) = w_{\hbar,0}(x,v) = \mathcal{F}_{y\rightarrow v}[\tilde{R}]_0(x,v),\label{eq:wigner equation data} 
\end{align}
with the pseudodifferential operator $\theta[V]$ defined as follows:
\begin{equation}\label{eq:def theta V}
    \theta[V]w_{\hbar}(x,v) := \tfrac{1}{(2\pi)^d}\int_{\mathbb{R}_y^d\times\mathbb{R}_{\eta}^d}e^{-i(v-\eta) \cdot y}\frac{V(x+\frac{\hbar}{2}y)-V(x-\frac{\hbar}{2}y)}{i\hbar}w_{\hbar}(t,x,\eta) \dd y \dd\eta.
\end{equation}
\begin{remark}Note that the structure of the symbol of $\theta[V]$ containing first order finite differences of the potential $V$ accounts for the convergence towards the classical force term of the Vlasov-type equation.
\end{remark}

\subsection{Velocity Averaging for the "quantum kinetic" Wigner equation}

For $\hbar$ fixed, B. Perthame and P.-L. Lions in 1992 \cite{lions1992lemmes} linked the dispersive structure of the free Schrödinger equation to moment lemmas for the (free) Wigner equation, viewed as a kinetic equation with $S[w_{\hbar}]\equiv 0$, i.e. the free transport equation for $w_{\hbar}$. The moment lemmas essentially express the fact that moments in $v$ of $w_{\hbar}$ combined with a certain decay in $x$ are more regular than $w_{\hbar}$, see for example Theorem 1 in \cite{lions1992lemmes}. In 1999, I. Gasser, P. Markowich and B. Perthame \cite{gasser1999dispersion} extended this result to the Schrödinger equation with potential by translating it to the Wigner equation with potential \eqref{eq:wigner equation} and applying the moment lemmas. This yields a $H^{1/2}$ gain in regularity for the Wigner function $w_{\hbar}$. While in the free case, the moment lemma of \cite{lions1992lemmes} can be translated into an averaging lemma by applying the moment lemma to $|\widehat{w_{\hbar}}|^2$ instead of $w_{\hbar}$, this correspondence is not obvious for the Wigner equation with potential. 

Hence in order to obtain an $L^2$-based averaging lemma for $w_{\hbar}$ we apply Theorem \ref{thm: diperna lions} to the Wigner equation directly. To this end we need an $L^2$ bound for $w_{\hbar}$, which is independent of $\hbar$ if and only if $R$ belongs a certain class of mixed states, a fact already observed in the semi-classical analysis of the Wigner-Poisson equation in \cite{lions1993mesures,markowich1993classical}. Indeed, by
Plancherel's theorem
\begin{equation}\label{WignerL2Bound}
\begin{aligned}
\iint|w_\hbar[R](x,v)|^2\dd x\dd v
=\frac1{(2\pi\hbar)^d}\iint|R(X,Y)|^2\dd X\dd Y= \frac1{(2\pi\hbar)^d}\sum_{j\ge 1}\lambda_j^2.
\end{aligned}
\end{equation}
Now assume that a family $\{R_\hbar\,:\,0<\hbar\le 1\}$ of density matrices satisfies
\[
\sup_{0<\hbar\le 1}\|w_\hbar[R]\|_{L^2}^2=C<\infty.
\]
By \eqref{WignerL2Bound}, this condition is implied by the following condition (stated in  in \cite{lions1993mesures,markowich1993classical})  on the occupation probabilities, 
\begin{equation}
    \sum_{j\geq 1} (\lambda_j^{\hbar})^2 \leq (2\pi \hbar)^d,
    \label{eq: condition probabilities}
\end{equation}
where we added a superscript $\hbar$ to emphasize that the eigenvalues $\lambda_j^{\hbar}$ depend on the semi-classical parameter $\hbar$. Condition \eqref{eq: condition probabilities}
implies that 
\[
\text{rank}R_\hbar\ge\frac1{(2\pi\hbar)^dC},\qquad\text{ so that }\varliminf_{\hbar\to 0}\left(\hbar^d\text{rank}R_\hbar\right)>0.
\]
In particular, this assumption rules out the possibility that $\{R_\hbar\,:\,0<\hbar\le 1\}$ is a family of pure states (i.e. rank-one density operators). Based on this assumption we can prove the following theorem on the "quantum averaging lemma" \eqref{eq:quantum lemma}, which only holds for $\hbar$ fixed. A "semi-classical averaging lemma", which holds uniformly in $\hbar$, will be stated in the next section. 

\begin{theorem}\label{thm:VAWignerL2 quantum}\textbf{\emph{(Quantum Averaging Lemma)}}
Let $V\in C(\mathbb{R}_x^d,\mathbb{R})$ such that $V\in L^\infty(\mathbb R^d_x,\mathbb{R})$ and let  $\{R_\hbar\,:\,\hbar\in(0,1]\}$ be a family of density operators $t\mapsto R_\hbar(t)$, $t\in[-T,T], T>0$, which are continuous on $[-T,T]$ for the weak operator topology, are weak solutions of the von Neumann equation \eqref{eq:von Neumann}
and satisfy the bound
\begin{equation}
	\sup_{|t|\le T}\tr_{L^2(\mathbb R^d)}(R_\hbar(t)^2)\le C^2(2\pi\hbar)^d
    \label{eq:bound hilbert schmidt}
\end{equation}
for some $C>0$. For each $\psi\in\mathcal S(\mathbb{R}_v^d)$, set
\[
	\rho_\psi[w_\hbar](t,x):=\int_{\mathbb R^d}w_\hbar(t,x,v)\psi(v)\dd v.
\]
Then for each $T>0$, there exists $C'_T>0$ such that
\begin{equation}
	\sup_{0<\hbar\le 1}\|\rho_\psi[w_\hbar]\|_{H^{1/2}((-T,T)\times\mathbb{R}_x^d)}\le C'_T\|V\|_{L^\infty(\mathbb{R}_x^d)}\hbar^{-1/2}.
    \label{eq:quantum lemma}
\end{equation}
\begin{proof}
   Write $\theta[V]w_{\hbar} = K[V] \ast_{v} w_{\hbar}$ where $$K[V] = \mathcal{F}_{y\rightarrow v} \left(\frac{V(x+\tfrac{\hbar y}{2})-V(x-\tfrac{\hbar y}{2})}{i\hbar}\right).\label{eq:convolution kernel}$$
    Then by the Plancherel theorem, $$\|\theta[V]w_{\hbar}\|_2 \leq \|\mathcal{F}[K[V]]\|_{\infty} \|w_{\hbar}\|_2.$$
    Then $\|\mathcal{F}[K[V]]\|_{\infty}$ can be bounded by $\tfrac{1}{\hbar} \|V\|_{L^{\infty}}$ times a constant, which proves the theorem by invoking Theorem \ref{thm: diperna lions} with $n=0$.
\end{proof}
\end{theorem}
Clearly, the constant on the RHS of \eqref{eq:quantum lemma} blows up as $\hbar\rightarrow 0$ and the estimate only holds for $\hbar$ fixed. Note that the regularity obtained, i.e. $\rho_{\psi}\in H^{1/2}$,  is the same as in the moment lemma from \cite{gasser1999dispersion}.  

\section{Semi-classical averaging lemmas}

\subsection{Semi-classical velocity averaging for mixed states}

We now present the following $L^2$-based "semi-classical averaging lemma" for mixed states satisfying condition \eqref{eq:bound hilbert schmidt}.

\begin{theorem}\label{thm:VAWignerL2 semi-classical}
 \textbf{\emph{(Semi-classical Averaging Lemma)}} Under the same assumptions as in Theorem \ref{thm:VAWignerL2 quantum}, in particular \eqref{eq:bound hilbert schmidt}, but with $V$ Lipschitz-continuous on $\mathbb R^d$ instead of $L^{\infty}$ we have that for each $T>0$,
\begin{equation}
	\sup_{0<\hbar\le 1}\|\rho_\psi[w_\hbar]\|_{H^{1/4}((-T,T)\times\mathbb{R}_x^d)}<\infty.
    \label{eq:semi-classical lemma}
\end{equation}

\end{theorem}
The proof of Theorem \ref{thm:VAWignerL2 semi-classical} involves more careful analysis of the convolution kernel $K[V]$, cf. \eqref{eq:convolution kernel}, in particular writing $K[V]$ as a divergence in $\xi$ (which leads to the exponent $1/4$ instead of $1/2$), 
but provides a uniform estimate in $\hbar$. We refer to \cite{golse2025velocity}  for a complete proof.

Contrary to Theorem \ref{thm:VAWignerL2 quantum}, the bound in \eqref{eq:semi-classical lemma} is uniform in (i.e. independent of) $\hbar$, thereby providing a uniform gain in regularity in the limit as $\hbar\rightarrow 0$. The price to pay is that this gain is 
weaker in the "semi-classical" case than in the "quantum" case ($H^{1/4}$ instead of $H^{1/2}$).


\subsection{Semi-classical velocity averaging for pure states}\label{sec:PureSt}


The semi-classical averaging lemma in \eqref{eq:semi-classical lemma} is a uniform in $\hbar$ estimate, based on uniform bounds in $L^2$ for the Wigner transform that are only available for a special class of mixed states satisfying 
\eqref{eq:bound hilbert schmidt}. Hence pure states are excluded by \eqref{eq:bound hilbert schmidt} and have to be dealt with separately.

As explained in Section \ref{sec:Prelim}, the density operator $R$ of a pure state 
can be written as
\begin{equation}\label{Rk1bra-ket}
	R = \ketbra{\psi}{\psi}, \quad R(X,Y) = \psi(X)\overline{\psi(Y)},
\end{equation}
where $R(X,Y)$ is the kernel of $R$, i.e. the density matrix.

As detailed in \cite{golse2025velocity}, pure states can be characterized in terms of a system of second order PDEs for the Wigner function. Denote by
\begin{equation}
    \tilde{w}_{\hbar}(x,y) := \mathcal F_{ v\to y}w_\hbar[R](x,y)
\end{equation}
the partial Fourier transform of 
$w_\hbar[R]$ in the $ v$-variable.
\begin{lemma}\label{L-CharRk1}
Let $R$ be a density matrix with Wigner function $w_\hbar[R](x, v)$ such that $\tilde{w}_{\hbar}\in C^1(\mathbb R^d\times\mathbb R^d)$ and assume that $\tilde{w}_{\hbar}\not=0$ for all $x,y\in\mathbb R^d$. 
Then $R$ is a rank-one density operator (i.e. a pure state) if and only if, for all $j,k=1,\ldots,d$,
\[
\left\{
\begin{aligned}
\frac4{\hbar^2}\partial_{y_j}\left(\frac{\partial_{y_k}\tilde{w}_{\hbar}}{\tilde{w}_{\hbar}}\right)
&=\partial_{x_j}\left(\frac{\partial_{x_k}\tilde{w}_{\hbar}}{\tilde{w}_{\hbar}}\right),
\\
\partial_{y_j}\left(\frac{\partial_{x_k}\tilde{w}_{\hbar}}{\tilde{w}_{\hbar}}\right)
&=\partial_{x_j}\left(\frac{\partial_{y_k}\tilde{w}_{\hbar}}{\tilde{w}_{\hbar}}\right),
\end{aligned}
\right.
\quad\text{ in }\mathcal D'(\mathbb{R}_x^d\times\mathbb{R}_y^d).
\]
\end{lemma}

Note that (a formal variant of) Lemma \ref{L-CharRk1} in the special case of one space dimension $d=1$ was already suggested by Tatarski\u\i\, in \cite{tatarskiui1983wigner}.

\begin{corollary}\label{C-KeyMadelungIdent}
Let $R=\ketbra{\psi}$, where $\psi\in C^1(\mathbb R^d)$ satisfies $\psi(x)\not=0$ for all $x\in\mathbb R^d$ and $\|\psi\|_{L^2(\mathbb R^d)}=1$. Let $\tilde{R}(x,y) = R(x+\tfrac{\hbar y}{2},x-\tfrac{\hbar y}{2})$.
Then
\[
\partial_{y_j}\partial_{y_k}\tilde{R}(x,y)=\frac{\partial_{y_j}\tilde{R}(x,y)\partial_{y_k}\tilde{R}(x,y)}{\tilde{R}(x,y)}
+\tfrac{\hbar^2}{4}\tilde{R}(x,y)\partial_{x_j}\left(\frac{\partial_{x_k} \tilde{R}(x,y)}{\tilde{R}(x,y)}\right)
\]
in the sense of distributions on $\mathbb R_x^d\times\mathbb R_y^d$, and for all $j,k=1,\ldots,d$.
\end{corollary}

\smallskip
Lemma \ref{L-CharRk1} and its consequence, Corollary \ref{C-KeyMadelungIdent} imply the following proposition, cf. \cite{golse2025velocity}. It gives sufficient conditions for which the Wigner measure $w$ of a pure state is monokinetic (and therefore fails to qualify for a semi-classical averaging lemma).

\begin{proposition}\label{P-Monokinetic}
Let $\psi_\hbar\in C^1(\mathbb R^d,\mathbb C\setminus\{0\})$ be a family of normalized wave functions satisfying the condition $\|\psi_\hbar\|_{L^2(\mathbb R^d)}=1$. Let  $w_\hbar:=w_\hbar[\ketbra{\psi_\hbar}]$.
Assume that 
\[
w_{\hbar_n}\to w\text{ in }\mathcal S'(\mathbb R^d\times\mathbb R^d),\quad\text{ and }\quad\hbar_n\|\nabla\rho_{\hbar_n}\|_{L^2(B(0,R))}\to 0
\]
for all $R>0$ as $\hbar_n\to 0$. Assume further that
\[
\mathcal E_{\hbar_n}=\tfrac1{2m}\int_{\mathbb{R}_v^d}|v|^2w_{\hbar_n}\dd v\to\mathcal E:=\tfrac1{2m}\int_{\mathbb{R}_v^d}| v|^2w\dd v\quad\text{ weakly in }L^2(\mathbb{R}_x^d)
\]
as $\hbar_n\to 0$, while
\[
\rho_{\hbar_n}=\int_{\mathbb{R}_v^d}w_{\hbar_n}\dd v\to\rho:=\int_{\mathbb{R}_v^d}w\dd v\quad\text{ and }J_{\hbar_n}=\tfrac1m\int_{\mathbb{R}_v^d} v w_{\hbar_n}\dd v\to J:=\tfrac1m\int_{\mathbb{R}_v^d} v w\dd v
\] 
strongly in $L^2_{loc}(\mathbb{R}_x^d)$ as $\hbar_n\to 0$ for each $R>0$. Set
\[
u(x):=\frac{\mathbf 1_{\rho(x)>0}}{\rho(x)}J(x)\in\mathbb{R}_x^d.
\]
Then $w$ is a monokinetic, positive Borel measure on the phase space $\mathbb{R}_x^d\times\mathbb{R}_v^d$, i.e.
\[
w=\rho(x)\delta( v-u(x)).
\]
\end{proposition}
That $w_{\hbar_n}$ converges to a monokinetic Wigner measure $w=\rho(x)\delta( v-u(x))$ is a very strong indication that the strong convergence of the first moments $\rho_{\hbar_n}$ and $J_{\hbar_n}$ in $L^2_{loc}(\mathbb{R}_x^d)$ cannot be deduced from a velocity averaging lemma. Indeed, velocity averaging is based on the fact that, at the kinetic level of description, the orthogonal projections of the variable
$\xi$ on each line through the origin are regularly distributed: see condition (2.1) in \cite{golse1988regularity}. This is incompatible with a situation where the $\xi$ variable concentrates on a single value at
each position $x$, which is precisely the case of a monokinetic distribution function.

\begin{remark}
By Theorem III.1 (5) in \cite{lions1993mesures}, the Wigner measure $\mu$ of a pure state is not necessarily monokinetic, which means that the additional assumptions in Proposition \ref{P-Monokinetic} are in general not satisfied by the family $\psi_\hbar$, and that these additional
assumptions are essential for Proposition \ref{P-Monokinetic} to hold. Moreover, the condition on $\nabla\rho_\hbar$ in Proposition \ref{P-Monokinetic} is obviously an essential feature of this result. We shall return to its physical meaning at the end of the next section. 
\end{remark}


\section{Connection to Quantum Hydrodynamics}\label{sec:DerMadelung}


Another application of Lemma \ref{L-CharRk1} and Corollary \ref{C-KeyMadelungIdent} leads to a  a quick and very natural derivation of the equations of Quantum Hydrodynamics (QHD). These equations
were obtained for the first time in \cite{Madelung1926} (see especially equations (4') and (3'') in that reference). In Madelung's own words, his work \cite{Madelung1926} corresponds mostly to an analogy
between Schr\"odinger's approach to quantum dynamics with hydrodynamics (see \cite{Madelung1926} on p. 322).

We assume that the wave function does not vanish at any point, i.e. $
	\psi(t,x)\not=0$  for all $x\in\mathbb R^d$ and $t>0$. The von Neumann equation 
    for $R(t)=\ketbra{\psi(t,\cdot)}$ is recast as
\begin{equation}
    \partial_t\tilde{R}(t,x,y)+\tfrac1m\sum_{j=1}^d \partial_{x_j}(-i\partial_{y_j})\tilde{R}(t,x,y)=\frac{V(x+\frac{\hbar}{2}y)-V(x-\frac{\hbar}{2}y)}{i\hbar}\tilde{R}(t,x,y),
    \label{eq:von Neumann Madelung}
\end{equation}
where $\tilde{R}(x,y) = R(x+\tfrac{\hbar y}{2},x-\tfrac{\hbar y}{2})$. Since $\psi$ is smooth and non-vanishing, the function $\tilde R(t,x,y)$ satisfies the assumptions of Lemma \ref{L-CharRk1}. 
We recall from \eqref{DefDensityFunc} that the density function is
\[
\rho(t,x):=\tilde{R}(t,x,0)=R(t,x,x)>0\qquad\text{ for all }x\in\mathbb R^d\text{ and }t>0,
\]
and define the velocity field $u(t,x)$ by the formula
\begin{equation}\label{DefVelField}
    u_j(t,x):=\frac{J_j(t,x)}{\rho(t,x)}=\frac{-i\partial_{y_j}\tilde R(t,x,0)}{m\tilde R(t,x,0)}, \quad j=1,\dots,d.
\end{equation}
Evaluating both sides of \eqref{eq:von Neumann Madelung} 
at $y=0$, we find  
the continuity equation
\begin{equation}
    \partial_t\rho(t,x)+\nabla_x\cdot(\rho(t,x)u(t,x))=0.
    \label{eq:continuity}
\end{equation}
For the Euler equation, we apply $-i\partial_{y_k}$ to both sides of \eqref{eq:von Neumann Madelung}, and arrive at the equality
\begin{equation}
    \partial_t(-i\partial_{y_k}\tilde{R}(t,x,y))-\tfrac1m\sum_{j=1}^d\partial_{x_j}(\partial_{y_j}\partial_{y_k} \tilde{R}(t,x,y))=-\partial_{y_k} \left(\tfrac{V(x+\frac{\hbar}{2}y)-V(x-\frac{\hbar}{2}y)}{\hbar}\tilde{R}(t,x,y)\right).
    \label{eq:euler 1}
\end{equation}
The second term on the left-hand side of \eqref{eq:euler 1} is transformed by using Corollary \ref{C-KeyMadelungIdent} so that
\[
\begin{aligned}
\tfrac1m\partial_t(-i\partial_{y_k}\tilde{R}(t,x,y))&+\tfrac1{m^2}\sum_{j=1}^d\partial_{x_j}\frac{(-i\partial_{y_j}\tilde{R}(t,x,y))(-i\partial_{y_k} \tilde{R}(t,x,y))}{\tilde{R}(t,x,y)} 
	\\
	&-\tfrac{\hbar^2}{4m^2}\sum_{j=1}^d\partial_{x_j}\left(\tilde{R}(t,x,y)\partial_{x_j}\left(\frac{\partial_{x_k}\tilde{R}(t,x,y)}{\tilde{R}(t,x,y)}\right)\right)
	\\
    	&+\tfrac1{2m}\left(\partial_kV\left(x+\tfrac{\hbar}{2}y\right)+\partial_kV\left(x-\tfrac{\hbar}{2}y\right)\right)\tilde{R}(t,x,y)
	\\
	&=\tfrac1m\tfrac{V(x+\frac{\hbar}{2}y)-V(x-\frac{\hbar}{2}y)}{\hbar}(-i\hbar\partial_{y_k}\tilde{R}(t,x,y)).
\end{aligned}
\]
Evaluating both sides of this equation at $y=0$ yields
\begin{equation}
\begin{aligned}
    \partial_t(\rho(t,x)u(t,x))&+\nabla_x\cdot(\rho(t,x)u(t,x)\otimes u(t,x))
    \\
    &=\tfrac{\hbar^2}{4m^2}\nabla_x\cdot(\rho(t,x)\nabla^2_x\ln\rho(t,x))-\tfrac1m\rho(t,x)\nabla V(x).
    \label{eq:Madelung Euler}
\end{aligned}
\end{equation}
The QHD equations in conservation form are \eqref{eq:continuity}-\eqref{eq:Madelung Euler}:
\begin{align}
\begin{cases}
&\partial_t\rho+\nabla_x\cdot(\rho u)=0,
	\\
	&\partial_t(\rho u)+\nabla_x\cdot(\rho u\otimes u)=-\tfrac1m\rho\nabla_xV + \tfrac{\hbar^2}{4m^2}\nabla_x\cdot(\rho\nabla^2_x\ln\rho).
\end{cases}
\label{eq:qhd}
\end{align}

There is another equivalent form of the Euler equation \eqref{eq:Madelung Euler}, 
\begin{equation}\label{eq:Madelung Euler QuantPress}
\partial_t(\rho u)+\nabla_x\cdot(\rho (u\otimes u+\tfrac1m\Pi))=-\tfrac1m\rho\nabla_xV
\end{equation}
with the \emph{quantum pressure tensor}
\begin{equation}\label{DefQuantPress}
\Pi:=-\tfrac{\hbar^2}{4m}\nabla^2\ln\rho, \quad (\nabla^2)_{jk} = \partial_{x_j} \partial_{x_k}
\end{equation}

Still another equivalent form of the QHD system requires checking the formula for the \emph{Bohm potential}. Observe that
\[
\begin{aligned}
-\tfrac{\hbar^2}{4m}\nabla_x\cdot(\rho\nabla^2_x\ln\rho) 
=:-\tfrac{\hbar^2}{2m}\rho\nabla_x\left(\frac{\Delta_x\sqrt\rho}{\sqrt\rho}\right) = \rho \nabla_x P.
\end{aligned}
\]
with the Bohm quantum potential 
\begin{equation}\label{BohmQuantPot}
P:=-\tfrac{\hbar^2}{2m}\frac{\Delta_x\sqrt\rho}{\sqrt\rho}\,,
\end{equation}
Then the Euler equation in the QHD system is recast as
\begin{equation}\label{eq:Madelung Euler BohmPot}
\partial_tu+u\cdot\nabla_xu=-\tfrac1m\nabla_x(V+P)\,.
\end{equation}
Since $\rho$ satisfies the continuity equation, equations \eqref{eq:qhd} and \eqref{eq:Madelung Euler QuantPress} are equivalent under the assumption that the density $\rho>0$ everywhere.

Finally, let us discuss the condition $\hbar^2\|\nabla\rho_\hbar\|^2_{L^2(\mathbb R^d)}\to 0$ in Proposition \ref{P-Monokinetic}. After straightforward calculations one obtains the identity 
\[
\hbar^2|\nabla\rho_\hbar|^2=\tfrac43\hbar^2\Delta(\rho^2_\hbar)+\tfrac83m\rho^2_\hbar P_\hbar=\hbar^2\Delta(\rho^2_\hbar)+2m\rho^2_\hbar\mathrm{Tr}\Pi_\hbar
\]

Owing to Proposition \ref{P-Monokinetic}, observe that $\hbar^2_n\rho^2_{\hbar_n}\to 0$ since $\rho_{\hbar_n}$ is bounded in $L^2(B(0,R))$. Hence
\[
\hbar^2_n\Delta(\rho^2_{\hbar_n})\to 0\quad\text{ in }\mathcal D'(\mathbb R^d)\qquad\text{ as }\hbar_n\to 0.
\]
Therefore
\[
\hbar_n^2\|\nabla\rho_{\hbar_n}\|^2_{L^2(B(0,R)}\to 0\iff\rho^2_{\hbar_n}P_{\hbar_n}\to 0\text{ in }\mathcal D'(\mathbb R^d)\iff\rho^2_{\hbar_n}\mathrm{Tr}\Pi_{\hbar_n}\to 0\text{ in }\mathcal D'(\mathbb R^d).
\]
In other words, the condition $\hbar_n^2\|\nabla\rho_{\hbar_n}\|^2_{L^2(B(0,R)}\to 0$ is equivalent to the fact that the quantum pressure $\Pi_{\hbar_n}$, or the Bohm potential $P_{\hbar_n}$, multiplied by 
$\rho_{\hbar_n}^2$ converges to $0$ in the sense of distributions on $\mathbb R^d$.

\section*{Acknowledgement}
Financial support from the Austrian Science Fund (FWF) is acknowledged, via the SFB project 10.55776/F65, as well as the Schrödinger grant 10.55776/J4840.\\
Also, F.G. and J.M. acknowledge the hospitality and support of the Wolfgang Pauli Institute Vienna.


\bibliographystyle{abbrv}
\bibliography{wignerveraging}

\end{document}